\documentclass[11pt]{amsart}

\newtheorem{theorem}{Theorem}[section]
\newtheorem{proposition}{Proposition}[section]
\newtheorem{lemma}{Lemma}[section]

\theoremstyle{definition}
\newtheorem{remark}{Remark}[section]

\numberwithin{equation}{section}

\newcommand{\Int}{{\mathbb Z}}

\newcommand{\Com}{{\mathbb C}}

\begin{document}

\title[On Maslankas's representation] {On Maslanka's representation for the
Riemann zeta function}

\author{Luis B\'{a}ez-Duarte}

\date{11 July 2003}
\email{lbaezd@cantv.net}


\begin{abstract}
A rigorous proof is given of the hypergeometric-like
representation of the Riemann zeta function discovered by K.
Maslanka.
\end{abstract}

\maketitle

\section{Introduction}
In \cite{maslanka1} and \cite{maslanka2} K. Maslanka introduced the following
representation of $\zeta(s)$ valid for all $s\in\Com$

\begin{equation}\label{pseries}
(s-1)\zeta(s)=\sum_{k=0}^\infty A_k P_k\left(\frac{s}{2}\right), 
\end{equation}
where the $A_k$ are given by

\begin{equation}\label{coeffs}
A_k=\sum_{j=0}^k (-1)^j {k \choose j} (2j+1)\zeta(2j+2),
\end{equation}

and the $P_k$ are the so-called \textit{Pochhammer polynomials} defined by

\begin{equation}\label{polynomials}
P_k(s):=\prod_{r=1}^k \left(1-\frac{s}{r}\right), \ \ \ (P_0(s):=1). 
\end{equation}

We apologize for having changed the notation used in \cite{maslanka1},
\cite{maslanka2}, but we have found it more natural to write it as above with future
applications in mind. 

Maslanka indicates that he was led to this interesting expression
by the desire to interpolate $\zeta(s)$ with a series of the form (\ref{pseries})
treating the $A_k$ as indeterminates to be found from the system of equations

$$
(2n-1)\zeta(2n)=\sum_{k=0}^\infty A_k P_k(2n), \ \ \ (n=1,2,\dots).
$$
This system is a triangular system because $P_k(m)=0$ when the integer $m\leq k$.
Its unique solution is found by a nice combinatorial argument given in Appendix A of
\cite{maslanka1}. 

The author of \cite{maslanka1}, \cite{maslanka2} gives two formal
proofs that the series (\ref{pseries}) with the coefficients as defined by
(\ref{coeffs}) represents $(s-1)\zeta(s)$. However no estimate of
the size of the $A_k$ is given that would imply the convergence of the series in any
region of the complex plane, and which would be needed to justify several of the formal
interchanges of limits and series in his proofs. Professor Maslanka himself kindly
indicated this gap to me in a personal communication. I thought it would be worthwhile
to devote some effort, which proved quite rewarding, to provide the missing steps in the
proof. An effective estimate of the rate at which $A_k\rightarrow 0$ is the crucial
missing link in the proof, and to this effect we shall show below that

\begin{equation}\label{estimate}
A_k\ll k^{-p}
\end{equation}
for every positive real. 

In closing this introduction we would like to point out an interesting paradox that 
could deserve some attention. Once the Maslanka representation is shown to be valid in
the whole complex plane it is obvious that (\ref{pseries}) will be valid for
$s=1,0,-1,-2,\dots$, noting however that in this case the series does not truncate 
one obtains interesting identities as pointed out in \cite{maslanka1}. If, on the 
other hand, one proposed at the outset the analogous problem of finding a
representation of the form

$$
(s-1)\zeta(s)=\sum_{k=0}^\infty c_k P_k(2-s),
$$
then it is clear that the series truncates for $s=1, 0, -1, -2,\dots$ giving again a
triangular system which leads to a surprising solution, namely

$$
(s-1)\zeta(s)=1+\frac{1}{2} (s-1)+\sum_{k=1}^\infty B_{k} P_k(2-s),
$$
where the $B_k$ are the Bernoulli numbers. This happens to be, in different garb, the
\textit{divergent} Euler-Maclaurin series of $\zeta(s)$. It may be objected that
Dirichlet series (albeit $(s-1)\zeta(s)$ is \textit{not} a Dirichlet series) are
determined by its values at any sequence with real parts  tending to $+\infty$, while
that is not the case toward the left. It may seldom be the case  that this question is
at all posed as ``all" Dirichlet series of interest have a finite abscissa  of
convergence so there is no way to move toward $-\infty$. On the other hand Professor
Malanska claims he has a representation that interpolates symultaneously at points
$2n$ and $3-2n$. We do not know what that series may be,  perhaps he means something
like $\sum c_k P_k(s/2)P_k((1-s)/2)$. 

These considerations shall be taken up in a more general study \cite{baez2} where 
we shall explore ``\textit{p-series}" $\sum c_k P_k(s)$ in general, and in particular in
\cite{baez1}  where we produce a new necessary and sufficient criterion for the Riemann
hypothesis (it looks like a discrete version of the well-known F. Riesz criterion
\cite{riesz}), namely, if

$$
b_k:=\sum_{j=0}^k (-1)^j {k \choose j} \frac{1}{\zeta(2j+2)}, \ \ \ (k=0, 1, 2, \dots),
$$ 
then the Riemann hypothesis is equivalent to the condition

$$
b_k\ll k^{-\frac{3}{4}+\epsilon}, \ \ \ (\forall \epsilon>0).
$$
This criterion is well-suited for numerical calculations. We have obtained a
very smooth and convincing smooth curve $b_k k^{\frac{4}{3}}$ tending to zero like
$\log^{-2}k$.

\section{Statement and proof of the main theorem}
We now formally state and prove the representation theorem.

\begin{theorem}\label{main}
For $A_k$ defined as in (\ref{coeffs}) we have
\begin{equation}\label{repr}
(s-1)\zeta(s)=\sum_{k=1}^\infty A_k P_k\left(\frac{s}{2}\right),
\end{equation}
for all $s\in\Com$. The convergence of the series is uniform and absolute in every
compact set of the complex plane.
\end{theorem}
Throughout it shall be very important to bear in mind the simple estimate for the
size of the Pochhammer polynomials contained in the following lemma. It is given here
without proof as it is quite standard.

\begin{lemma}\label{estpoch}
For every disk $|s|\leq r <\infty$ there is a positive constant $C_r$ such that
\begin{equation}\label{pochmajor}
|P_k(s)|\leq C_r k^{-\Re s}, \ \ \ (|s|\leq R).
\end{equation}
\end{lemma}

Do note also the equally trivial facts contained in:

\begin{lemma}\label{pochandgamma}
\begin{equation}\label{pochgamma}
P_k(s)=(-1)^k{s-1\choose k}=\frac{\Gamma(k+1-s)}{k!\Gamma(1-s)}.
\end{equation}
\end{lemma}
\ \\

\begin{proof}[Proof of Theorem \ref{main}]
The key estimate (\ref{estimate}) for the coefficients $A_k$ shall be proved in
the next section. We accept it for the moment and show now how the informal steps in
Maslanka's second proof can be made rigorous. We shall first show that with no
assumption on the size of the $A_k$ one can prove equality in (\ref{pseries}) for
$\Re s>4$. We start with Maslanka's clever idea of expressing $s-1$ as the
derivative of $-\alpha^{1-s}$ at $\alpha=1$. We indicate by $D_\alpha$ the operation of
differentiating with respect to $\alpha$ and proceed as follows:
 
\begin{eqnarray}\label{trick}\nonumber
\alpha^{-s}(s-1)\zeta(s) &=&
-D_\alpha\alpha^{1-s}\sum_{n=1}^\infty\frac{1}{n^s}\\\nonumber
&=&-D_\alpha\alpha\sum_{n=1}^\infty\frac{1}{(\alpha n)^s}\\\nonumber
&=&-D_\alpha\alpha\sum_{n=1}^\infty\frac{1}{(\alpha n)^2}
\left(\frac{1}{(\alpha n)^2}\right)^{\frac{s}{2}-1}\\\nonumber
&=&-D_\alpha\alpha\sum_{n=1}^\infty\frac{1}{(\alpha
n)^2}\left(1-\left(1-\frac{1}{(\alpha n)^2}\right)\right)^{\frac{s}{2}-1}\\\nonumber
&=&-D_\alpha\alpha\sum_{n=1}^\infty\frac{1}{(\alpha n)^2}
\sum_{k=0}^\infty (-1)^k {\frac{s}{2}-1 \choose k}
\left(1-\frac{1}{(\alpha n)^2}\right)^k\\
&=&-D_\alpha\sum_{n=1}^\infty\sum_{k=0}^\infty
P_k\left(\frac{s}{2}\right)\frac{1}{\alpha\hspace{.5mm} n^2} 
\left(1-\frac{1}{(\alpha n)^2}\right)^k,
\end{eqnarray}
where we have applied the binomial theorem assuming as we may that the derivative is
calculated from the right so $\alpha > 1$. Consider now the formally differentiated
double series on the right-hand side above, namely

$$
\sum_{n=1}^\infty\sum_{k=0}^\infty
P_k\left(\frac{s}{2}\right) 
\frac{1}{(\alpha n)^2}\left(1-\frac{2k+1}{(\alpha n)^2}\right)
\left(1-\frac{1}{(\alpha n)^2}\right)^{k-1}.
$$
Using Lemma \ref{estpoch} we see that this double series is uniformly termwise
majorized by a convergent double series of positive terms in $\Re s \geq 4+2\epsilon$
since

$$
\left|P_k\left(\frac{s}{2}\right)\right|\frac{1}{\alpha\hspace{.5mm} n^2}
\left(1-\frac{2k+1}{(\alpha n)^2}\right) 
\left(1-\frac{1}{(\alpha\hspace{.5mm} n)^2}\right)^k
\leq k^{-1-\epsilon}\frac{1}{n^2}. 
$$
Thus we can both invert the sums at the end of (\ref{trick}) and do the
differentiation termwise at $\alpha=1$ obtaining

\begin{equation}
(s-1)\zeta(s)=\sum_{k=0}^\infty P_k\left(\frac{s}{2}\right)
\sum_{n=1}^\infty
\frac{1}{n^2}\left(1-\frac{2k+1}{n^2}\right)
\left(1-\frac{1}{n^2}\right)^{k-1}.
\end{equation}
We easily calculate the last series on the right as

$$
\sum_{j=0}^{k-1}(-1)^j{k-1\choose j}(\zeta(2j+2)-(2k+1)\zeta(2j+4))
$$
which, after separating in two sums and shifting indices in the second, becomes

$$
\sum_{j=0}^{k}(-1)^j{k-1\choose j}(2j+1)\zeta(2j+2),
$$
which is none other than $A_k$. One could even argue this would be a direct way to
find the coefficients. We have thus shown the desired relation (\ref{pseries}) at least
for $\Re s > 4$. However the strong estimate (\ref{estimate}), to be proved below,
providing that $A_k\ll k^{-p}$ for every $p>0$, together with Lemma \ref{estpoch}
immediately shows that the series $\sum A_k P_k(s/2)$ converges uniformly on any
compact of the plane, thus defining an entire function that must be equal to
$(s-1)\zeta(s)$ by analytic continuation.  
\end{proof}

\section{The size of the coefficients.}
We shall establish here this Proposition:

\begin{proposition}\label{aestimate}
For any $p>0$ there is a constant $C_p>0$ such that

\begin{equation}\label{Asmall}
|A_k|\leq C_p k^{-p}.
\end{equation}
\end{proposition}

To lighten up the proof we first prove some lemmas. Define the crucial sequence of
rational functions $\phi_k(x)$ by

\begin{equation}\label{phik}
\phi_k(x):=\left(1-\frac{1}{x^2}\right)^k\frac{1}{x}, \ \ \ (k\in\Int^+).
\end{equation}
We now record a very simple while extremely important property of $\phi(x)$.
\begin{lemma}\label{phivanish}
For $0\leq a \leq k$ denote
$$
\phi^{(a)}_k(x):=\frac{d^a}{dx^a}\phi_k(x),
$$
then we have
\begin{equation}\label{0der}
\phi^{(a)}_k(1)=\phi^{(a)}_k(\infty)=0,
\end{equation}
where the second equality is unrestricted on $a$. All the $\phi^{(a)}_k(x)$ for
$a\geq1$ are integrable.
\end{lemma}
The following lemma is the essential and quite elementary tool to get a grip on
the size of $A_k$.

\begin{lemma}\label{dershape}
Define a double sequence of polynomials $p_{a,j}$, with integral coefficients as
follows: Let $p_{0,0}=1$, $p_{a,j}=0$ when $j<0$ and when $j>a$, and specify the
recurrence equation

\begin{equation}\label{precc}
p_{a,j}=-(2j+a)p_{a-1,j}+(2k+2j-a)p_{a-1,j-1}.
\end{equation}
With these definitions we have for any integers $k\geq a \geq 0$

\begin{equation}\label{aderexpr}
\phi^{(a)}_k(x)=\left(1-\frac{1}{x^2}\right)^{k-a}\sum_{j=0}^a \frac{p_{a,j}(k)}{x^{a+2j+1}}.
\end{equation}
\end{lemma}
The simple proof is by (double) induction. The recurrence relation (\ref{precc}) will not be used
other than the fact that it is part of the proof. It is nice to record it explicitly for further
use if need be.

Having straightened some notation, as all the above was essentially, we are ready
for a first preliminary estimate.

\begin{lemma}\label{normder}
For any fixed integer $a\geq1$ and $\epsilon>0$ there is a constant $C=C(a,\epsilon)$
such that

\begin{equation}\label{estnormder}
\int_1^\infty |\phi_k^{(a)}(x)|dx\leq C k^{-\frac{a}{2}(1-\epsilon)}.
\end{equation}

\end{lemma}

\begin{proof}
Take $k\geq a\geq1$. We split the interval of integration at the point
$x=k^{\frac{1}{2}-\frac{\epsilon}{3}}$. We express the $a$-th derivative according to Lemma
\ref{dershape} equation (\ref{aderexpr}). For the finite range of the integral we have

\begin{eqnarray}\label{range1}\nonumber
\int_1^{k^{\frac{1}{2}-\frac{\epsilon}{3}}}|\phi^{(a)}_k(x)|dx
&\leq&\int_1^{k^{\frac{1}{2}-\frac{\epsilon}{3}}}\left|\left(1-\frac{1}{x^2}\right)^{k-a}
\sum_{j=0}^a\frac{p_{a,j}(k)}{x^{a+2j+1}}\right|dx\\\nonumber
&\leq&\left|\sum_{j=1}^{a}p_{a,j}(k)\right|
\int_1^{k^{\frac{1}{2}-\frac{\epsilon}{3}}}\left(1-\frac{1}{x^2}\right)^{k-a}dx\\\nonumber
&\leq&
O_a(k^a)\left(1-\frac{1}{k^{1-\frac{2\epsilon}{3}}}
\right)^{k-a}k^{\frac{1}{2}-\frac{\epsilon}{3}}\\
&=&O_a(e^{-k^{\frac{\epsilon}{3}}}).
\end{eqnarray}
Now for the infinite range of integration we have

\begin{eqnarray}\label{range2}\nonumber
\int_{k^{\frac{1}{2}-\frac{\epsilon}{3}}}^\infty|\phi^{(a)}_k(x)|dx
&\leq&\int_{k^{\frac{1}{2}-\frac{\epsilon}{3}}}^\infty
\left|\left(1-\frac{1}{x^2}\right)^{k-a}
\sum_{j=0}^a\frac{p_{a,j}(k)}{x^{a+2j+1}}\right|dx\\\nonumber
&\leq&
\sum_{j=0}^a |p_{a,j}(k)|\int_{k^{\frac{1}{2}-\frac{\epsilon}{3}}}^\infty
\frac{1}{x^{a+2j+1}}dx\\\nonumber
&\leq&
\sum_{j=0}^a
\frac{1}{a+2j}\frac{|p_{a,j}(k)|}{k^{(\frac{1}{2}-\frac{\epsilon}{3})(a+2j)}}
\\\nonumber
&\leq&
\frac{1}{k^{a(\frac{1}{2}-\frac{\epsilon}{3})}}\sum_{j=0}^a
\frac{1}{a+2j}\frac{|p_{a,j}(k)|}{k^{j(1-\frac{2\epsilon}{3})}}\\\nonumber
&\leq&
\frac{1}{k^{a(\frac{1}{2}-\frac{\epsilon}{3})}}\sum_{j=0}^a
\frac{1}{a+2j}O_a(k^{\frac{2\epsilon}{3}j})\\
&\leq&
\frac{1}{k^{a(\frac{1}{2}-\frac{\epsilon}{3})}}O_a(k^{\frac{2\epsilon}{3}a})
=O_a\left(k^{-\frac{a}{2}(1-\frac{\epsilon}{2})}\right).
\end{eqnarray}
Since this last estimate (\ref{range2}) surely dominates the smaller one
in (\ref{range1}) we can allow the substitution of $\epsilon/2$ by $\epsilon$ to arrive
at the desired (\ref{estnormder}).    
\end{proof}
\begin{remark}\label{sufforder}
Actually there is no need for $\epsilon$ to be arbitraly small for the use
this estimate (\ref{estnormder}) is destined for. If one takes, say $\epsilon=1$ we
see that the order obtained is $O_a(k^{-\frac{a}{4}})$, which is plenty, since we
are planning to differentiate a large but fixed number of times, while
$k\rightarrow\infty$. The idea perhaps was to see how much we can eke out of this
method, which is as close as $k^{-\frac{a}{2}}$ as desired. Perhaps a more painstaking
analysis, taking into account the actual nature of the polynomials $p_{a,j}$, could
yield an essentially faster order of convergence to zero, which is desirable for the
effectivenes of the Maslanka representation.
\end{remark}
At this point it is high time to connect the coefficients $A_k$ with the function
$\phi_k$. So we have first a high school triviality, namely:

\begin{lemma}\label{binsum1}
For any $x\not=0$
\begin{equation}\label{binomeq}
\sum_{j=0}^k (-1)^j {k \choose j}(2j+1)\frac{1}{x^{2j+2}}=-\phi^{'}_k (x).
\end{equation}
\end{lemma}
\begin{proof}
\begin{eqnarray}\nonumber
\sum_{j=0}^k (-1)^j {k \choose j}(2j+1)\frac{1}{x^{2j+2}}&=&
-\frac{d}{dx}
\left(\frac{1}{x}\sum_{j=0}^k (-1)^j {k \choose j}\frac{1}{x^{2j}}\right)\\\nonumber
&=&-\frac{d}{dx}
\left(\frac{1}{x}\left(1-\frac{1}{x^2}\right)^k\right).
\end{eqnarray}
\end{proof}

It should now be clear why the work on the derivatives of $\phi_k(x)$ if the reader
notices that the sum below should be  estimated by Euler-Maclaurin with a
large number of terms.

\begin{lemma}\label{Aasphi}
\begin{equation}\label{muntz}
A_k=-\sum_{n=1}^\infty \phi^{'}_k (n).
\end{equation}
\end{lemma}

\begin{proof}
The following interchange of sums is totally elementary to justify:

\begin{eqnarray}\nonumber
A_k &=& \sum_{j=0}^k (-1)^j {k \choose j}(2j+1)\zeta(2j+2)\\\nonumber
&=&\sum_{j=0}^k (-1)^j {k \choose j}(2j+1)\sum_{n=1}^\infty \frac{1}{n^{2j+2}}\\\nonumber
&=&\sum_{n=1}^\infty\sum_{j=0}^k (-1)^j {k \choose j}(2j+1)\frac{1}{n^{2j+2}}\\\nonumber
&=&-\sum_{n=1}^\infty \phi^{'}_k (n),
\end{eqnarray}
where, of course, we applied Lemma \ref{binsum1} in the last equality.
\end{proof}

We are finally ready to prove the estimate for the $A_k$.

\begin{proof}[Proof of Proposition \ref{aestimate}]
Take an integer $a>4p$ and assume that $k>a$; naturally one will make
$k\rightarrow\infty$. The crucial  elementary properties of $\phi_k$ and its
derivatives expressed in Lemma \ref{phik} imply the remarkable fact that the
application of the Euler-Maclaurin summation formula to a depth of $a$ steps to the sum
in (\ref{muntz}), that is to

$$
A_k=-\sum_{n\geq1}\phi^{'}_k(n)
$$
results in the sum being equal to the remainder term! Therefore

$$
A_k=-\frac{(-1)^a}{a!}\int_1^\infty \overline{B_a}(x)\phi^{(a)}_k(x)dx,
$$
where $\overline{B_a}(x)$ is the $a$-th periodified Bernoulli polynomial. Now
apply Lemma \ref{normder} eq. (\ref{estnormder}) with $\epsilon=\frac{1}{2}$ to get

\begin{eqnarray}\nonumber
|A_k|&\leq&\frac{\|\overline{B_a}\|_\infty}{a!}\int_1^\infty|\phi_k^{(a)}(x)|dx\\\nonumber
&\leq&O_a(k^{-\frac{a}{4}})=O_a(k^{-p}).
\end{eqnarray}

\end{proof}

\begin{remark}
We believe there is room for improvement in estimating the smallness of the $A_k$.
This should have a bearing on how useful the Maslanka formula may turn out to be in
the end.
\end{remark}

\bibliographystyle{amsplain}
 
 \ \\
\ \\
\noindent Luis B\'{a}ez-Duarte\\
Departamento de Matem\'{a}ticas\\
Instituto Venezolano de Investigaciones Cient\'{\i}ficas\\
Apartado 21827, Caracas 1020-A\\
Venezuela\\
\email{lbaezd@cantv.net}

\end{document}